\documentclass[13pt]{article}
\usepackage{latexsym}

\usepackage{authblk}

\usepackage[dvips]{graphicx}
\usepackage[normalem]{ulem}
\usepackage{psfrag,amsmath,bm,psfrag,pstricks,xspace}
\usepackage[a4paper,top=1cm,bottom=1.1cm,left=1cm,right=1cm]{geometry}

\usepackage[francais,english]{babel}
\usepackage{amssymb,amsthm,epsfig,amsmath,color,hyperref,url,xspace}

\begin{document}

\newcommand{\refe}[1]{Ref.~\cite{#1}}
\newtheorem{teo}{Theorem}[section]
\newenvironment{preuve}{\textbf{\noindent Proof: \\*}\rm}{$\hfill\square$ \vspace{4mm}}
\newtheorem{defi}[teo]{Definition}
\newtheorem{corol}[teo]{Corollary}
\newtheorem{prop}[teo]{Proposition}
\newtheorem{theorem}[teo]{Theorem}
\newtheorem{criterion}[teo]{Criterion}  
\newtheorem{lemma}[teo]{Lemma} 
\newtheorem{examp}[teo]{Example}

\newcommand{\varpoly}{X}
\newcommand{\cofvo}{v_0}
\newcommand{\cofv}{v}
\newcommand{\eqcardan}{\left( \mathcal{C}\right) }
\newcommand{\zcardan}{z}
\newcommand{\pcardan}{p}
\newcommand{\qcardan}{q}
\newcommand{\deltacardan}{\Delta_{\mathcal{C}}}
\newcommand{\indcardan}{n}
\newcommand{\zindcardan}{\zcardan_{\indcardan}}
\newcommand{\rootcardan}{R_{\mathcal{C}}}
\newcommand{\prefaca}{\alpha_1}
\newcommand{\prefacb}{\alpha_2}
\newcommand{\eqxm}{x_{m}}

\newcommand{\matm}{\mathbf{M}} 
\newcommand{\matma}{\mathbf{M}_1} 
\newcommand{\matmb}{\mathbf{M}_2} 

\newcommand{\vecva}{\mathbf{v}_0}
\newcommand{\vecv}{\mathbf{v}}

\newcommand{\compma}{m_1} 
\newcommand{\compmb}{m_2} 
\newcommand{\compmc}{m_3} 

\newcommand{\compmd}{m_4} 
\newcommand{\compme}{m_5} 
\newcommand{\compmf}{m_6} 

\newcommand{\compmg}{m_7} 
\newcommand{\compmh}{m_8} 
\newcommand{\compmi}{m_9} 

\newcommand{\compya}{y_0} 
\newcommand{\compza}{z_0}
\newcommand{\compx}{v_1} 
\newcommand{\compy}{v_2} 
\newcommand{\compz}{v_3}

\newcommand{\vala}{\lambda_1} 
\newcommand{\valb}{\lambda_2} 
\newcommand{\valc}{\lambda_3} 
\newcommand{\valp}{k} 
\newcommand{\valk}{\tilde{k}} 
\newcommand{\valka}{k_0} 

\newcommand{\para}{\tilde{a}}
\newcommand{\parb}{\tilde{b}}
\newcommand{\parm}{m_0}
\newcommand{\paro}{o}
\newcommand{\parl}{l}

\newcommand{\coa}{A}
\newcommand{\cob}{B}
\newcommand{\coc}{C}
\newcommand{\cod}{D}
\newcommand{\coe}{E}

\newcommand{\coda}{d_1}
\newcommand{\codb}{d_2}
\newcommand{\codc}{d_3}
\newcommand{\codd}{d_4}
 
\newcommand{\cofa}{f_1}
\newcommand{\cofb}{f_2}
\newcommand{\cofc}{f_3}
\newcommand{\cofd}{f_4}

\newcommand{\cora}{r_1}
\newcommand{\corb}{r_2}
\newcommand{\corc}{r_3}
\newcommand{\cord}{r_4}
\newcommand{\core}{r_5}
\newcommand{\corf}{r_6}

\newcommand{\delp}{\Delta_{\valp}}
\newcommand{\cofea}{e_1}
\newcommand{\cofeb}{e_2}
\newcommand{\cofec}{e_3}
\newcommand{\cofed}{e_4}
\newcommand{\cofca}{c_1}
\newcommand{\cofcb}{c_2}
\newcommand{\cofna}{s_1}
\newcommand{\cofnb}{s_4}
\newcommand{\cofnc}{s_5}
\newcommand{\cofnd}{s_7}
\newcommand{\cofne}{s_8}
\newcommand{\cofnf}{s_9}
\newcommand{\cofng}{s_{10}}
\newcommand{\cofnh}{s_6}
\newcommand{\cofni}{s_2}
\newcommand{\cofnj}{s_3}

\newcommand{\cofnind}{n_i}
\newcommand{\ind}{i}
\newcommand{\cofla}{l_1}
\newcommand{\coflb}{l_2}
\newcommand{\coflc}{l_3}
\newcommand{\cofld}{l_4}
\newcommand{\cofle}{l_5}
\newcommand{\coflf}{l_6}
\newcommand{\coflg}{l_7}
\newcommand{\coflh}{l_8}
\newcommand{\cofli}{l_9}
\newcommand{\colb}{C_{\parl}}
\newcommand{\colc}{B_{\parl}}
\newcommand{\cold}{D_{\parl}}

\newcommand{\ctels}{s_{\parl}}
\newcommand{\rl}{r_{m}}
\newcommand{\ensol}{\mathcal{S}_{\parl}}

\newcommand{\tr}{\mathrm{Tr}}

\newcommand{\esc}{\mathbb{C}}  

\newcommand{\eqsym}{\left( \mathcal{S}\right) } 

\newcommand{\eqpoly}{\left( \mathcal{P}\right) } 
\newcommand{\eqx}{x}
\newcommand{\eqb}{b}
\newcommand{\eqc}{c}
\newcommand{\eqd}{d} 
\newcommand{\eqs}{s}

\newcommand{\mata}{\mathbf{A}}
\newcommand{\delo}{\Delta_{\paro}} 
\newcommand{\dell}{\Delta_{\parl}} 
\newcommand{\denumo}{d_{\paro}}
\newcommand{\lo}{l_{\paro}}

\newcommand{\eq}[1]{Eq.~\eqref{#1}}
\newcommand{\Eq}[1]{Equation~\eqref{#1}}
\newcommand{\ie}{i.e.\@\xspace}

\newcommand{\valpa}{\valp_1}
\newcommand{\valpb}{\valp_2}

\newcommand{\parba}{\parb\left(\valpa \right) }
\newcommand{\parbb}{\parb\left(\valpb \right) }

\newcommand{\parla}{\parl\left(\valpa \right)}
\newcommand{\parlb}{\parl\left(\valpb \right)}

\newcommand{\indma}{m_1}
\newcommand{\indmb}{m_2}
\newcommand{\indm}{m}

\newcommand{\eqxc}{x_0}
\newcommand{\eqcc}{c_0}
\newcommand{\eqdc}{d_0}
\newcommand{\deloc}{\delta_{\paro}}
\newcommand{\dellc}{\delta_{\parl}}

\newcommand{\ac}{\mathcal{A}}
\newcommand{\bc}{\mathcal{B}}
\newcommand{\cc}{\mathcal{C}}
\newcommand{\dc}{\mathcal{D}}
\newcommand{\mc}{\mathcal{M}}
\newcommand{\nc}{\mathcal{N}}

\newcommand{\arga}{A_1}
\newcommand{\argb}{A_2}
\newcommand{\raca}{\mathcal{R}_1}
\newcommand{\racb}{\mathcal{R}_2}
\newcommand{\deltsl}{\delta_{\parl}}

\title{Analytical formula for the roots of the general complex cubic polynomial }
\author[1]{Ibrahim Baydoun}
\date{\today}

\vspace{0.2cm}

\affil[1]{ESPCI ParisTech, PSL Research University, CNRS, Univ Paris Diderot, Sorbonne Paris Cit\'{e}, Institut Langevin, 1 rue Jussieu, F-75005, Paris, France}

\vspace{0.5cm}

\maketitle
\abstract{
We present a new method to calculate analytically the roots of the general complex polynomial of degree three. 
This method is based on the approach of appropriated changes of variable involving an arbitrary parameter. 
The advantage of this method is to calculate the roots of the cubic polynomial as uniform formula using the standard convention of the square and cubic roots. 
In contrast, the reference methods for this problem, as Cardan-Tartaglia and Lagrange, give the roots of the cubic polynomial as expressions with case distinctions which are incorrect using the standar convention. 
}
\\

\noindent\textbf{keyword}:~cubic polynomial roots, appropriated change of variable, analytical uniform formula.




  
\section{Introduction}

The problem of solution formulas of polynomial equations is fundamental in algebra \cite{VDW,SML-GB}. It is useful, by example, for the eigenvalues perturbation \cite{Ha-Pe} and the solution of a system of nonlinear equations \cite{IP-NA}. 
Moreover, it has many physical applications, as the study of the singularities of the surfaces of refractive indices in crystal optics \cite{K-M-S} and the diagonal eigenvalues perturbation \cite{And}. 
Another example is the calculation of the sound velocity anisotropy in cubic crystals \cite{Ts-Yo}, or the diagonalization of Christoffel tensor to calculate the velocities of the three quasi-modes of the elastic waves in anisotropic medium \cite{Bay}. 
An explicit calculation for the particular case of $3 \times 3$ real symmetric matrices has been studied in \cite{Me-Ra}. 

For a cubic equation, Cardano-Tartaglia and Lagrange's formulas are subject to different interpretations depending on the choices for the values of the square and cubic roots of the formulas \cite{Bour,Lag}. The correct interpretation is generally derived from conditions between the roots and the coefficients of the cubic polynomial. 
It is a fundamental question to be able to define uniform formula that yields correct solutions for all cubic equations. Ting Zhao's thesis and related work recently brought an important advance on Lagrange's formulas with uniform solution formulas for cubic and quartic equations with real coefficients. 
The results was obtained by introducing a new convention for the argument of the cubic root of a complex number \cite{Zh-Wa}. 
For cubic equations, we present in this paper a uniform formula that is valid with complex coefficients. This more general result are uniquely defined by the standard convention of the square and cubic roots. 

To illustrate the problem, we introduce Cardan-Tartaglia formula and afterward we give two examples. Indeed, Cardan-Tartaglia calculate the 3 solutions of the equation
\begin{align*}
\eqcardan: \zcardan^3 + \pcardan\zcardan + \qcardan = 0
\end{align*}
as 
\begin{align*}
\zindcardan := \exp\left( \indcardan\frac{2\pi}{3}\sqrt{-1}\right) \sqrt[3]{\frac{1}{2}\left(-\qcardan+\sqrt{\frac{\deltacardan}{27}} \right) }
         + \exp\left(- \indcardan\frac{2\pi}{3}\sqrt{-1}\right) \sqrt[3]{\frac{1}{2}\left(-\qcardan-\sqrt{\frac{\deltacardan}{27}} \right) };
        ~ \indcardan \in\left\lbrace 0,1,2\right\rbrace  ,
\end{align*}
where $\deltacardan$ is defined as 
\begin{align*}
\deltacardan : = -\left(  4 \pcardan^3 + 27 \qcardan^2\right) .
\end{align*}
We show by the following examples that the standard convention of the cubic roots
\begin{align*}
 \arg \left( \sqrt[3]{\varpoly} \right) = \frac{1}{3} \arg\left(  \varpoly\right) ; ~ \varpoly\in \esc
\end{align*}
is not relevant in the above Cardan-Tartaglia's expression. Indeed:
\begin{examp}\label{exam-1}
The solution of $\eqcardan$ for the case $\pcardan = 1$ and $\qcardan = 2$ is 
\begin{align*}
\zindcardan\frac{-1-\sqrt{-3}}{2}   ;
        ~ \indcardan \in\left\lbrace 0,1,2\right\rbrace  ,
\end{align*}
instead of $\zindcardan$.
\end{examp}
\noindent Example \ref{exam-1} requires to shift the standard convention by $-2\pi/3$ in order to be suitable for the Cardan-Tartaglia's solutions. Another example requires another convention:
\begin{examp}\label{exam-2} 
The solution of $\eqcardan$ for the case $\pcardan = 1$ and $\qcardan = \sqrt{-1}$ is 
\begin{align*}
\zindcardan \frac{-1+\sqrt{-3}}{2}  ;
        ~ \indcardan \in\left\lbrace 0,1,2\right\rbrace  ,
\end{align*}
instead of $\zindcardan$.
\end{examp}
\noindent Therefore, a natural question is how we can calculate the roots of $\eqcardan$ using uniquely the standard convention? 
\newline
To overcome this problem, a mathematical software as Matlab generates the symbolic solutions of $\eqcardan$ in terms of one cubic root $\rootcardan$ as follows 
\begin{multline*}
\rootcardan -\frac{\pcardan}{3\rootcardan} , 
\\
\frac{\sqrt{-3}}{2} \left( \frac{\pcardan}{3\rootcardan} + \rootcardan \right)
+ \frac{\pcardan}{6\rootcardan} - \frac{\rootcardan}{2} , 
\\
  - \frac{\sqrt{-3}}{2}\left( \frac{\pcardan}{3\rootcardan} + \rootcardan \right) + \frac{\pcardan}{6\rootcardan}
- \frac{\rootcardan}{2},
\end{multline*}
where $\rootcardan$ is given by
\begin{align*}
\rootcardan = \sqrt[3]{ \sqrt{\frac{-\deltacardan}{2^2 3^3}}  - \frac{\qcardan}{2} }.
\end{align*}
On the one hand, this formulation hid the problem of the convention of the cubic roots. But, on the other hand, it involves another problem by dividing by $\rootcardan$ which leads to a singularity in the formulation if $\rootcardan$ is zero.
So, the above question can be reworded by the following manner: how we calculate the roots of $\eqcardan$ using uniquely the standard convention and without any singularity? 
\newline
The response to this question is the objective of this paper.

We present a uniform analytical formula for the roots of the general complex cubic polynomial.
Our method is based on the approach of changes of variable involving an arbitrary parameter, which rests on two principal ideas to diagonalize the general $3 \times 3$ complex matrix: 
\\
Firstly, we introduce one change of variable involving an arbitrary parameter. So, the characteristic equation of this matrix can be replaced by another non-polynomial equation in terms of the first new variable.
\\
Secondly, we introduce another change of variable involving another arbitrary parameter. The second variable replaces the first and the second arbitrary parameter can be chosen as required to make the form of the new equation as the sum of a cube of a certain monomial and a term independent of the unknown second variable. 
\\
Therefore, we can solve easily this new polynomial equation and consequently deduce the solution of the original equation, which is the eigenvalues of the general $3 \times 3$ complex matrix.

In the sequel, we apply this approach to calculate the roots of the general complex cubic polynomial. In fact, we construct a particular matrix such that its characteristic polynomial will coincide with this cubic polynomial. Consequently, the eigenvalues of this particular matrix become identical to the roots of this cubic polynomial.

The rest of the paper is structured as follows. In section \ref{sec-1}, we derive the approach of appropriated changes of variable involving an arbitrary parameter. Consequently in section \ref{sec-2}, we apply it on the general complex cubic polynomial.

\section{Approach of appropriated changes of variable involving an arbitrary parameter}\label{sec-1}

Let $\matm$ be the general $3 \times 3$ complex matrix:
\begin{align*}
\matm:=\begin{pmatrix}
   \compma & \compmb & \compmc \\
   \compmd & \compme & \compmf \\
   \compmg & \compmh & \compmi 
\end{pmatrix}.
\end{align*}
In this section, we aim to calculate the spectrum of $\matm$ with this approach. We detail it by the following four steps. 
The first step consists of performing a change of variable for the eigenvalues of $\matm$ so that we get a non-polynomial equation. 
In the second step, we will deal with this equation in order to simplify it by some algebraic manipulations.
The third step will introduce a second change of variable in order to replace the non-polynomial equation by another polynomial equation.
Finally, in the fourth step, we will solve the latter polynomial equation to deduce the spectrum of $\matm$.

Starting the first step by introducing the following definition which will be useful in the sequel:
\begin{defi}\label{defe-coef-eq-car}
We associate to $\matm$ the following expressions in terms of its components:
\begin{align*}
&\cofca:=\compma-\compme-\frac{\compmc\compmh}{\compmb} ,  \hspace{10mm} 
\cofcb:=\cofca^2+4\left( \compmb\compmd+\compmc\compmg\right) ,  \nonumber\\
&\cofea:=\frac{3}{2}\left(\compma-\cofca \right) -\frac{1 }{2}\tr\left(\matm \right) ,  \hspace{10mm}
\cofeb:=\frac{3 }{2}\left(\compma-\cofca \right)^2 +\frac{3}{8}\left(\cofcb-\cofca \right)^2 - \left(\compma-\cofca \right)\tr\left(\matm \right) +\frac{1 }{4}\left\{\left[\tr\left( \matm\right)\right] ^2-\tr\left( \matm^2\right) \right\} ,  \nonumber\\
&\cofec:=\frac{1}{8}\left[\left(2\compma-\cofca \right)^3+3\left(2\compma-\cofca \right)\cofcb \right]  -\frac{1}{4} \left[ \left(2\compma-\cofca \right)^2+\cofcb\right] \tr\left(\matm \right) +\frac{1 }{4} \left(2\compma-\cofca \right) \left\{\left[ \tr\left( \matm\right) \right] ^2 -\tr\left( \matm^2\right)\right\}-\det\left( \matm\right)  ,  \nonumber\\
&\cofed:=\frac{ 1 }{8}\left[ 3\left(2\compma-\cofca \right)^2+\cofcb \right]   -\frac{1 }{2} \left(2\compma -\cofca \right) \tr\left(\matm \right)  +\frac{1 }{4}\left\{\left[\tr\left( \matm\right)\right] ^2-\tr\left( \matm^2\right)\right\}.
\end{align*}
$\tr\left(\matm \right) $ and $\det\left(\matm \right)$ stand respectively for the trace and the determinant of matrix $\matm$.
\end{defi}
\noindent Now, we introduce the first change of variable involving an arbitrary parameter so that the eigenvalues of $\matm$ can be rewritten in terms of this new variable $\parb$ such that $\parb$ verifies a non-polynomial equation like the following proposition:
\begin{prop}\label{p-a-b-prop}
An eigenvalue $\valp$ of matrix $\matm$ can be written in terms of two unknowns $\para$ and $\parb$ as follows:
\begin{align}\label{p-a-b}
\valp = \frac{\para-\parb}{2} + \frac{1}{2}\left[2\compma - \cofca  + \sqrt{\delp } \right] ,
\end{align}
where $\para$ can be chosen arbitrary and $\parb$ is governed by the following non-polynomial equation:
\begin{align}\label{p-equ-a-b}
\coa+\cob\parb+\coc\parb^2+\cod\sqrt{\delp}+\coe\parb\sqrt{\delp} +\frac{1}{2}\left(\parb^2\sqrt{\delp}-\parb^3 \right) =0,
\end{align}
where $\delp$ is defined by
\begin{align} \label{p-equ-a-b-cof}
\delp := \parb^2+2\left(\cofca-\para \right) \parb +\para^2-2\cofca\para+\cofcb = 
\left(\para-\parb - \compma+\compme +\frac{\compmc\compmh}{ \compmb} \right)^2 +4\left( \compmb\compmd+\compmc\compmg\right) .
\end{align}
The coefficients of \eq{p-equ-a-b} depend on parameter $\para$ as follows:
\begin{multline} \label{cof-c-e}
\coa = \frac{\para^3}{2}+\cofea\para^2+\cofeb\para+\cofec ,    \hspace{5mm}                  
\cob = -\frac{3\para^2}{2}-2\cofea\para-\cofeb ,               \hspace{5mm}                  
\coc = \frac{3\para}{2}+\cofea ,                                                \\
\cod  = \frac{\para^2}{2}+\left( \cofea+\frac{\cofca}{2}\right) \para+\cofed ,  \hspace{5mm} 
\coe  = -\para-\left( \cofea +\frac{\cofca}{2}\right) ,
\end{multline}
where the coefficients of \eq{cof-c-e} are given by definition \ref{defe-coef-eq-car}.
\end{prop}
\begin{preuve}
Firstly, we decompose the matrix $\matm$ as the sum of two matrix $\matma$ and $\matmb$ as follows: 
\begin{align}\label{decomp-matm}
\matm=\underbrace{\begin{pmatrix}
   \para & 0     & 0      \\
   0     & \parb & -\parm \\
   0     & 0     & 0 
\end{pmatrix}}_{\matma}
+\underbrace{\begin{pmatrix}
   \compma -\para & \compmb        & \compmc \\
   \compmd        & \compme -\parb & \compmf +\parm \\
   \compmg        & \compmh        & \compmi 
\end{pmatrix}}_{\matmb},
\end{align}
where $\parm$ is chosen such that $\matmb$ has an eigenvector of the form $[0,\compya,\compza]$. The advantage of this choice is to construct an eigenvector of $\matm$ in terms of these of $\matmb$ which will be elaborated in the sequel.
Indeed, if we take $\parm$ such that
\begin{align*}
\parm = \left( \compme-\parb + \compmc \frac{\compmh}{\compmb}  -\compmi \right) \frac{\compmc}{\compmb} -\compmf ,
\end{align*}
then $\matmb$ have the eigenvalue  
\begin{align}\label{ka-a-b}
\valka = -\compmc \frac{\compmh}{\compmb}  +\compmi
\end{align}
relative to the eigenvector 
\begin{align*}
\vecva:=\left[ 0,-\frac{\compmc}{\compmb },1\right] .
\end{align*}
Equation \ref{ka-a-b} gives one eigenvalue $\valka$ of $\matmb$, then the characteristic polynomial of $\matmb$
\begin{align*}
\det\left[ \varpoly \begin{pmatrix}
   1 & 0 & 0 \\
   0 & 1 & 0 \\
   0 & 0 & 1
\end{pmatrix} - \matmb \right] 
\end{align*}
is divided in $\esc\left[ \varpoly\right] $ by the polynomial $\varpoly - \valka$.
Consequently, we deduce that the two other eigenvalues of $\matmb$ verify a second order polynomial and they are given by: 
\begin{align}\label{k-a-b}
\valk =  \frac{-\para-\parb}{2} + \frac{1}{2}\left[2\compma - \cofca\mp\sqrt{ \delp } \right]  .
\end{align}
Secondly, we search an eigenvector for $\matm$ of the form
\begin{align*}
\vecva + \cofv\vecv;~ \cofv \in \esc,
\end{align*}
where $\hat{\vecv}:=[\compx,\compy,\compz]$ is an unit eigenvector of $\matmb$ relative to $\valk$. 
Indeed, using 
\begin{align*}
& \matmb\vecva = \valka \vecva , \\
& \matmb\vecv = \valk \vecv , 
\end{align*}
we obtain from \eq{decomp-matm} that the vectorial equation
\begin{align*}
\matm\left( \vecva + \cofv\vecv \right) =\valp \left( \vecva + \cofv\vecv \right);~ \valp \in \esc, 
\end{align*}
is equivalent to the following system:
\begin{displaymath}
\eqsym \left\{ \begin{array}{ll}
 \compx \cofv \left( \para+\valk - \valp \right) = 0  , \\
\parb\left(\cofv \compy - \frac{\compmc}{\compmb } \right) - \parm\left( \cofv\compz+1\right) +\cofv\valk\compy - \frac{\compmc}{\compmb }\valka  - \valp\left( \cofv \compy - \frac{\compmc}{\compmb } \right) =0 ,  & \textrm{} \\
\cofv\valk\compz + \valka - \valp\left(\cofv \compz + 1 \right) =0  . & \textrm{}
\end{array} \right.
\end{displaymath}
The impact of the choice of $\parm$ appears in the fact that the first equation of $\eqsym$ is independent of $\vecva + \cofv\vecv$ and it can be verified by taking $\valp = \para+\valk$. Thus, by taking $\cofv$ as 
\begin{align*}
\cofv  = \frac{1}{\compz}\frac{\valp - \valka}{ \valk - \valp   },
\end{align*}
the third equation of $\eqsym$ is verified. It remains to verify the second equation of $\eqsym$. Indeed, $\valp=\para+\valk$ and \eq{k-a-b} implies \eq{p-a-b} where $\valp$ depends only on $\para-\parb$ since $\delp$ depends only on $\para-\parb$.
Therefore, by inserting \eq{p-a-b} in the characteristic equation of $\matm$, we get
\begin{align}\label{p-equ}
\valp ^3-\tr\left(\matm \right) \valp ^2+\frac{1 }{2}\{\left[ \tr\left( \matm\right)\right] ^2-\tr\left( \matm^2\right)\}\valp - \det(\matm)=0, 
\end{align}
so that we can choose arbitrarily among $\para$ and $\parb$ one parameter as required, while the other parameter will be the unique unknown of \eq{p-equ}, as well as $\valp$ is an eigenvalue of $\matm$. 
We obtain that if $\para$ and $\parb$ satisfy \eq{p-equ}, then the second equation of $\eqsym$ is verified. 
To end the proof, we develop \eq{p-equ} to deduce that \eq{p-equ} and \eq{p-equ-a-b} are equivalent.
\end{preuve}

In the second step, we aim to simplify \eq{p-equ-a-b} of the unknown $\parb$, where $\para$ can be chosen as required. So, we derive from \eq{p-equ-a-b} two equations of $\parb$ given by the following two propositions:
\begin{prop}\label{d-equ-a-b-prop}
For fixed $\para$, $\parb$ satisfies the following equation:
\begin{align}\label{d-equ-a-b}
\coda+\codb\parb+\codc\parb^2+\codd\sqrt{\delp}-\codc \parb\sqrt{\delp}=0,
\end{align}
where the coefficients of \eq{d-equ-a-b} are given by:
\begin{align*} 
\coda:=\codc\para^2+\cofna\para+\cofni,     \hspace{5mm}                 
\codb:=-2\codc\para-\cofna,                 \hspace{5mm}                     
\codc:=\compmb\compmd+\compmc\compmg,       \hspace{5mm}
\codd:=\codc\para+\cofnj.
\end{align*}
Here $\cofna$, $\cofni$ and $\cofnj$ are given using definition \ref{defe-coef-eq-car} as follows:
\begin{align}\label{na-equ-cof}
 \cofna := \frac{\cofca\cofcb}{2} +\cofca\cofeb+\cofcb\cofea-2\cofca\cofed-\cofec, &\hspace{5mm} 
&\cofni := \cofca\cofec+\cofcb\cofed , &\hspace{5mm}
&\cofnj := \cofec+\cofca\cofed.
\end{align}
\end{prop}
\begin{preuve}
\eq{d-equ-a-b} is deduced by multiplying \eq{p-equ-a-b} by $\parb + \sqrt{\delp} - \para +\cofca$.
\end{preuve}
\begin{prop}\label{f-equ-a-b-prop}
For fixed $\para$, $\parb$ satisfies the following equation:
\begin{align}\label{f-equ-a-b}
\cofa+\cofb\parb+\cofc\parb^2+\cofd\sqrt{\delp}+\cofc \parb\sqrt{\delp}=0,
\end{align}
where the coefficients of \eq{f-equ-a-b} are given by:
\begin{align*} 
&\cofa:=\cofc\para^2+\cofnb\para+\cofnc,     & \hspace{5mm}                 
&\cofb:=-2\cofc\para-\cofnb,                                    \nonumber\\
&\cofc:=\left(\compme-\compma \right) \compmc\compmg + \left(\compmi-\compma \right) \compmb\compmd 
-\compmb\compmf \compmg -\compmc \compmd\compmh  ,   & \hspace{5mm}
&\cofd:=-\cofc\para + \cofnh.
\end{align*}
Here $\cofnb$, $\cofnc$ and $\cofnh$ are given using definition \ref{defe-coef-eq-car} as follows:
\begin{align}\label{nb-equ-cof}
&\cofnb := \cofcb\codc -2\cofca^2\cofed -2\cofca\cofec -2 \codc\cofeb, & \hspace{5mm} 
&\cofnc := \cofcb\cofec-2\codc\cofec +\cofca\cofcb\cofed,              & \hspace{5mm} 
&\cofnh := \cofca\cofec +\cofcb\cofed -2\codc\cofed.
\end{align}
\end{prop}
\begin{preuve}
Similary to the proof of proposition \ref{d-equ-a-b-prop}, on one hand we multiply \eq{p-equ-a-b} by $2\codc$ and on the other hand we multiply \eq{d-equ-a-b} by $\sqrt{\delp}$. Consequently, we subtract the obtained equations to deduce \eq{f-equ-a-b}.
\end{preuve}
\\
Now, we show that equations (\ref{d-equ-a-b}) and (\ref{f-equ-a-b}) imply the two equations (\ref{r-equ-a-b-bis}) and (\ref{r-equ-a-b}) given by the following corollary:
\begin{corol}
For fixed $\para$, $\parb$ satisfies the following equations:
\begin{align}
& \cora+\corb\parb+\corc\sqrt{\delp} -2 \parb\sqrt{\delp}=0, \label{r-equ-a-b-bis} \\
& \cord+\core\parb+\corf\sqrt{\delp} +2 \parb^2=0 ,          \label{r-equ-a-b}
\end{align}
where the coefficients of \eq{r-equ-a-b-bis} and \eq{r-equ-a-b} are given by:
\begin{align}\label{r-equ-cof}
&\cora:=-\corb\para +\cofnf ,                           & \hspace{5mm}                 
&\corb:=-\frac{\cofna}{\codc}+\frac{\cofnb}{\cofc},     & \hspace{5mm}                                
&\corc:=2\para+\cofng ,  \nonumber\\
&\cord:= 2\para^2+\cofnd \para+\cofne, & \hspace{5mm}
&\core:= -4\para-\cofnd  , & \hspace{5mm}
&\corf:= \frac{\cofnj}{\codc}+\frac{\cofnh}{\cofc}
\end{align}
Here, the coefficients of \eq{r-equ-cof} are defined using definition \ref{defe-coef-eq-car}, \eq{na-equ-cof} and \eq{nb-equ-cof} as follows:
\begin{align*} 
& \cofnd := \frac{\cofna}{\codc}+\frac{\cofnb}{\codc},  & \hspace{5mm} 
&\cofne := \frac{\cofni}{\codc}+\frac{\cofnc}{\cofc},   & \hspace{5mm} 
&\cofnf := \frac{\cofni}{\codc}-\frac{\cofnc}{\cofc} ,  & \hspace{5mm}
&\cofng:=\frac{\cofnj}{\codc}-\frac{\cofni-2\codc\cofed}{\cofc} .
\end{align*}
\end{corol}
\begin{preuve}
We divide \eq{d-equ-a-b} and \eq{f-equ-a-b} respectively by $\codc$ and $\cofc$. Then, the subtraction and the addition of the obtained equations imply \eq{r-equ-a-b-bis} and \eq{r-equ-a-b}.
\end{preuve}

In the third step, we introduce the second change of variable
\begin{align}\label{chang-var-l-o}
\parl:=-\left(\parb+\paro\sqrt{\delp} \right) 
\end{align}
in order to deduce from the intractable equations of $\parb$ another polynomial quation of the new variable $\parl$, where $\paro$ is an arbitrary parameter can be chosen as required:
\begin{prop}\label{l-equ-o-prop}
Let $\paro\in\esc$. Then $\parl$ verifies the following equation:
\begin{align}\label{l-equ-o}
2\parl^3+\colc\parl^2+\colb\parl+\cold=0,
\end{align}
where the coefficients of \eq{l-equ-o} are given by:
\begin{align*} 
&\colc:=-\cofla\paro+6\para-\coflb ,
&\colb:=\coflc\paro^2-2\left(\cofla\para+\cofld \right) \paro+6\para^2-2\coflb\para+\cofle,
 \nonumber \\
&\cold:=-\coflf\paro^3+\left(\coflc\para-\coflg \right)\paro^2-\left(\cofla\para^2+2\cofld\para+\coflh \right)\paro 
+2\para^3-\coflb\para^2+\cofle\para+\cofli  .
\end{align*}
Here $\colc$, $\colb$ and $\cold$ are defined using definition \ref{defe-coef-eq-car}, \eq{na-equ-cof}, \eq{nb-equ-cof} and \eq{r-equ-cof} as follows:
\begin{align*} 
&\cofla:=2\frac{\cofna}{\codc}-\cofnd+\corf,            \hspace{10mm}
\coflb:=\corf-\cofnd-2\frac{\cofnj}{\codc},             \hspace{10mm}
\coflc:=\cofne-2\cofcb+2\cofca\corf+\frac{\corf\cofna-\left( \cofnd+4\cofca\right) \cofnj }{\codc},
\nonumber \\
&\cofld:=2\frac{\cofni}{\codc}-\cofne-\cofca\corf,                        \hspace{10mm}
\cofle:= \cofne+\frac{\cofnd \cofnj-\corf\cofna}{\codc},                  \hspace{10mm}
\coflf:=2\cofca\cofne+\cofcb\cofnd+\frac{\cofna\left(\cofne-2\cofcb \right)-\left(4\cofca+\cofnd\right) \cofni }{\codc},
\nonumber \\
&\coflg:=2\cofca\cofne+\cofcb\left( \cofnd+\corf\right) +\frac{\left( \cofne-2\cofcb\right)\cofnj -\corf\cofni  }{\codc},                                                                 \hspace{10mm}
\coflh:=\cofcb\corf+\frac{\cofnd\cofni -\cofna\cofne}{\codc},             \hspace{10mm}
\cofli:=\frac{\cofne \cofnj -\corf\cofni }{\codc}.
\end{align*}
\end{prop}
\begin{preuve}
The first part of the proof consists in getting expressions for $\parb^2$ and $\sqrt{\delp}$ just in terms of $\parl$ and $\paro$. Indeed, \eq{chang-var-l-o} implies:
\begin{align}\label{pre-b-ol}
\parb = -\left( \paro\sqrt{\delp}+\parl\right).
\end{align}
We take the square of \eq{pre-b-ol} and we replace $\delp$ by its value given in \eq{p-equ-a-b-cof}. Then, we insert \eq{pre-b-ol} in the obtained equation to deduce an expression of $\parb^2$. This chain of operations is illustrated as follows:
\begin{multline}\label{pre-b2-ol}
\parb^2 =  \paro^2 \underbrace{\delp}_{\parb^2+2\left(\cofca-\para \right) \parb +\para^2-2\cofca\para+\cofcb} + \parl^2 + 2\paro\parl\sqrt{\delp} \Rightarrow 
\left(1- \paro^2\right) \parb^2 =  2\paro^2 \left(\cofca-\para \right) \underbrace{\parb}_{ -\left( \paro\sqrt{\delp}+\parl\right)} + \paro^2\left(  \para^2-2\cofca\para+\cofcb \right)  + \parl^2 + 2\paro\parl\sqrt{\delp} 
 \\
\Rightarrow  \parb^2 = \frac{ 2\paro \left[\parl- \left(\cofca-\para \right) \paro^2 \right]\sqrt{\delp}+\left(\para^2-2\cofca\para+\cofcb\right)\paro^2+\parl^2-2\left(\cofca-\para \right)\paro^2\parl}{1-\paro^2} .
\end{multline}
So, inserting \eq{pre-b-ol} and \eq{pre-b2-ol} in \eq{r-equ-a-b} in order to substitute $\parb$ and $\parb^2$ and consequently to get the expression of $\sqrt{\delp}$ as follows:
\begin{multline}\label{pre-del-ol}
\cord+\core\underbrace{\parb}_{-\left( \paro\sqrt{\delp}+\parl\right)}+\corf\sqrt{\delp} 
+2 \underbrace{\parb^2}_{ \left(1-\paro^2 \right) ^{-1}\left\lbrace  2\paro \left[\parl- \left(\cofca-\para \right) \paro^2 \right]\sqrt{\delp}+\left(\para^2-2\cofca\para+\cofcb\right)\paro^2+\parl^2-2\left(\cofca-\para \right)\paro^2\parl \right\rbrace }=0 
 \\
\Rightarrow \sqrt{\delp}=\frac{\left(1-\paro^2 \right)\left(\cord-\core\parl \right)+2\left[ \left( \para^2-2\cofca\para+\cofcb\right) \paro^2-2\left(\cofca-\para \right) \paro^2\parl+\parl^2\right]   }
{\left(1-\paro^2 \right)\left(\core\paro-\corf \right)+2\left[ 2\left(\cofca-\para \right) \paro^3-2\paro\parl\right]}.
\end{multline}
The second part of the proof consists to use the previous expressions of $\parb^2$ and $\sqrt{\delp}$ in order to deduce \eq{l-equ-o}. Indeed, we take the square of the equation $\parb+\paro\sqrt{\delp}+\parl=0$ deduced from \eq{chang-var-l-o}. Then, we insert \eq{p-equ-a-b-cof} and \eq{d-equ-a-b} in the obtained equation to substitute $\delp$ and $\parb\sqrt{\delp}$. Afterwards we insert \eq{pre-b-ol} to substitute $\parb$. This chain of operations is illustrated as follows:
\begin{multline*}
\parb^2+\paro^2\underbrace{\delp}_{\parb^2+2\left(\cofca-\para \right) \parb +\para^2-2\cofca\para+\cofcb }+\parl^2 + 2\paro \underbrace{\parb\sqrt{\delp}}_{ \codc^{-1}\left( \coda+\codb\parb+\codc\parb^2+\codd\sqrt{\delp}\right) } + 2\parl\paro\sqrt{\delp} + 2 \parl\parb = 0 
\Rightarrow \\
\left( 1+\paro\right)^2 \parb^2  
+ 2 \left[  \left(\cofca-\para \right)\paro^2  +  \parl +\frac{\codb}{\codc}\paro \right] \underbrace{\parb}_{-\left( \paro\sqrt{\delp}+\parl\right)}
+ 2\left( \parl + \frac{\codd}{\codc} \right) \paro \sqrt{\delp}
+ \left( \para^2-2\cofca\para+\cofcb \right) \paro^2  +\parl^2 + 2\frac{\coda }{\codc} \paro = 0
\Rightarrow \\
\left( 1+\paro\right)^2 \parb^2 
+ 2 \left[   \frac{\codd}{\codc}  - \left(\cofca-\para \right)\paro^2  -\frac{\codb}{\codc}\paro \right]\paro\sqrt{\delp}
- 2 \left[  \left(\cofca-\para \right)\paro^2  +  \parl +\frac{\codb}{\codc}\paro \right]\parl
+ \left( \para^2-2\cofca\para+\cofcb \right) \paro^2  +\parl^2 + 2\frac{\coda }{\codc} \paro = 0 .
\end{multline*}
Finally, to complete the proof, it is enough to insert \eq{pre-b2-ol} and \eq{pre-del-ol} in the above equation in order to substitute $\parb^2$ and $\sqrt{\delp}$ and consequently to deduce \eq{l-equ-o}.
\end{preuve}
\\
Since $\paro$ is an arbitrary parameter, then we can choose $\paro$ as required to solve equation (\ref{l-equ-o}) of the unknown $\parl$: 
\begin{lemma}\label{lem-1}
If $\paro$ verifies the following polynomial equation of degree two:
\begin{align}\label{equ-o-teo}
\left(\cofla^2-6\coflc \right) \paro^2+2\left( \cofla\coflb+6\cofld\right) \paro +\coflb^2-6\cofle=0, 
\end{align}
then $\parl=-\para+\rl$, where
\begin{align*} 
\rl = \frac{\cofla\paro+\coflb}{6} +\indm \ctels ; ~~~\indm \in \left\{ -1 ~,~~  \frac{1-\sqrt{-3}}{2} ~,~~ \frac{1+\sqrt{-3}}{2}\right\}.
\end{align*}
Here $\ctels$ is given by:
\begin{align*} 
\ctels:=\sqrt[3]{ \left(\frac{\cofla\paro+\coflb}{6} \right)^3-\frac{1}{2} \left(\coflf\paro^3+\coflg\paro^2+\coflh\paro-\cofli \right) }. 
\end{align*}
\end{lemma}
\begin{preuve}
Equation (\ref{l-equ-o}) can be rewritten as follows
\begin{align}\label{l-equ-teo-o}
\parl^3+3\frac{\colc}{6}\parl^2+3\frac{\colb}{6}\parl+\left(\frac{\colc}{6} \right)^3 +\ctels^3
-\frac{\left(\cofla^2-6\coflc \right) \paro^2+2\left( \cofla\coflb+6\cofld\right) \paro +\coflb^2-6\cofle}{12}\para=0.
\end{align}
The condition $\colc^2=6\colb$, which is equivalent to \eq{equ-o-teo}, allows us to factorise \eq{l-equ-teo-o} as follows:
\begin{align}\label{l-equ-teo-o-fac}
\left( \parl+\frac{\colc}{6}\right)^3 +\ctels^3= \left( \parl+\frac{\colc}{6}+\ctels\right) \left[\left( \parl+\frac{\colc}{6}\right)^2-\ctels\left( \parl+\frac{\colc}{6}\right)+ \ctels^2\right] =0.
\end{align}
So, \eq{l-equ-teo-o-fac} can be solved easily, wich ends the proofs. 
\end{preuve}
\\
Lemma \ref{lem-1} gives the expressions for the set of solutions of $\parl$. Then, we deduce from the change of variable (\ref{pre-b-ol}) the following polynomial equation of $\parb$, which will be useful to give an explicit formula for the set of solutions of $\parb$:
\begin{corol}
If $\paro$ verifies \eq{equ-o-teo}, then $\parb$ verifies the following equation
\begin{align}\label{equ-cor1-b}
 \left(1-\paro^2 \right)\parb^2+2\left[ \rl-\para+\paro^2\left(\para-\cofca \right)\right]\parb+\left(\rl-\para \right) ^2-\paro^2 \left(\para^2-2\cofca\para+\cofcb  \right) =0 ,
\end{align}
where
\begin{align*} 
\rl \in \ensol:= \left\{ \frac{\cofla\paro+\coflb}{6}-\ctels ~;~~  \frac{\cofla\paro+\coflb}{6}+\frac{1-\sqrt{-3}}{2}\ctels ~;~~ \frac{\cofla\paro+\coflb}{6}+\frac{1+\sqrt{-3}}{2}\ctels\right\}.
\end{align*}
\end{corol}
\begin{preuve}
We have from \eq{chang-var-l-o} that
\begin{align}\label{b-del-cor1}
\sqrt{\delp}=-\frac{1}{\paro}\left(\parb +\parl\right). 
 \end{align}
The proof ends by taking the square of \eq{b-del-cor1} and afterwards replacing $\delp$ by its value given by \eq{p-equ-a-b-cof} in the obtained equation, where $\parl$ is replacing by $\rl-\para$ from lemma \ref{lem-1}.
\end{preuve}

Finally, in the fourth step, we can determine the set of solutions of $\parb$ by solving the nonlinear system constructed by \eq{equ-cor1-b} and one of the following equations (\ref{d-equ-a-b}), (\ref{f-equ-a-b}), (\ref{r-equ-a-b-bis}) or (\ref{r-equ-a-b}) of $\parb$. Thus, we can solve this nonlinear system by choosing $\para$ as required to simplify the calculation. Consequently, we deduce from \eq{p-a-b} the spectrum of $\matm$:
\begin{corol}\label{cor2-p}
The spectrum of matrix $\matm$ is given by:
\begin{align}\label{sol-cor2-p}
\valp=\frac{1}{2}\left[\left( 1+\frac{1}{\paro} \right) \frac{\codc\rl^2+\cofnj\left( 1-\paro\right)\rl-\cofni\paro\left( 1-\paro\right)-\codc \cofcb \paro^2 }
{\left( 1-\paro\right)\left( \cofnj+\cofna\paro\right) -2\codc\cofca\paro^2+ \codc\left( 1+\paro\right)\rl} +2\compma-\cofca-\frac{\rl}{\paro}\right] ,
\end{align}
where $\rl\in\ensol$ and $\paro$ verifies \eq{equ-o-teo}.
\end{corol}
\begin{preuve}
Inserting \eq{b-del-cor1} in \eq{d-equ-a-b}, we get an equation of $\parb$ and $\parb^2$. Then, we substitute $\parb^2$ using  \eq{equ-cor1-b} to obtain:
\begin{align}\label{b-cor2}
\parb=\para- \frac{\codc\rl^2+\left( 1-\paro\right)\cofnj\rl-\cofni\paro\left( 1-\paro\right)-\codc \cofcb \paro^2 }
{\left( 1-\paro\right)\left( \cofnj+\cofna\paro\right) -2\codc\cofca\paro^2+ \codc\left( 1+\paro\right)\rl}.
\end{align}
Finally, to deduce \eq{sol-cor2-p}, we insert \eq{b-del-cor1} in \eq{p-a-b} in order to substitute $\sqrt{\delp}$ and afterwards we insert \eq{b-cor2} in the obtained equation in order to substitute $\parb$. 
\end{preuve}
\\
The denumerator of the fraction in \eq{sol-cor2-p} contains $\rl$ which involves a cubic root. So, we calculate another expression for spectrum of $\matm$ in order to rationalize this denumerator. The two expressions of spectrum of $\matm$ will be useful in the proof of theorem \ref{teo-centrale}.
\begin{corol}\label{cor3-p}
The spectrum of matrix $\matm$ is given by:
\begin{align}\label{sol-cor3-p}
\valp=\frac{1}{2}\left[\left( 1+\frac{1}{\paro} \right) \frac{\cofnf\paro^2-\cofne\paro+\left( \corf-\cofng\paro\right)\rl  } {\left(\cofnd-\cofng\right)\paro +\corb\paro^2+\corf+2\paro\rl  }  +2\compma-\cofca-\frac{\rl}{\paro}\right] ,
\end{align}
where $\rl\in\ensol$ and $\paro$ verifies \eq{equ-o-teo}.
\end{corol}
\begin{preuve}
We insert \eq{b-del-cor1} in both \eq{d-equ-a-b} and \eq{f-equ-a-b} to replace $\sqrt{\delp}$ by its value in terms of $\parb$, $\parl$ and $\paro$. Then, by substituting $\parb^2$ from one obtained equation in the other, we deduce that:
\begin{align}\label{b-cor3}
\parb=\para- \frac{\cofnf\paro^2-\cofne\paro+\left( \corf-\cofng\paro\right)\rl  } {\left(\cofnd-\cofng\right)\paro +\corb\paro^2+\corf+2\paro\rl  }.
\end{align}
Finally, to deduce \eq{sol-cor3-p}, we insert \eq{b-del-cor1} in \eq{p-a-b} to substitute $\sqrt{\delp}$ and afterwards we insert \eq{b-cor3} in the obtained equation to substitute $\parb$. 
\end{preuve}

\section{Application on the cubic polynomial}\label{sec-2}

Let $\eqpoly$ be the general complex polynomial of degree three in $\esc\left[ \eqx \right] $:
\begin{align*} 
\eqpoly: \eqx^3+\eqb\eqx^2+\eqc\eqx+\eqd=0;~~  \eqb,\eqc,\eqd\in\esc.
\end{align*}

We aim to calculate analytically the roots of $\eqpoly$ by applying the results of section \ref{sec-1}. So, we introduce the following definitions which will be useful to find the expressions for the roots of $\eqpoly$:
\begin{defi}\label{defe-del-equ}
We associate to $\eqpoly$ the following expressions in terms of its coefficients:
\begin{align*} 
&\dell:= 2\eqc^3\left( 8\eqb^6+132\eqb^3\eqd+36\eqd^2+\eqc^3+33\eqb^2\eqc^2-66\eqb\eqc\eqd\right)  
+12\eqb^4\eqc \left( \eqd^2 -7\eqc^3\right) 
 -\eqb^2 \eqc^2 \eqd\left( 24\eqb^3 +291\eqd\right)
 +\eqd^3\left( 144\eqb\eqc-2\eqb^3 -27\eqd\right) , \nonumber \\ 
&\delo:= -4\eqb^3\eqd+\eqb^2\eqc^2+18\eqb\eqc\eqd-4\eqc^3-27\eqd^2 , 
\end{align*}
and
\begin{align*} 
\denumo:=4\eqb^4\eqc^2-4\eqb^3\eqc\eqd-14\eqb^2\eqc^3+\eqb^2\eqd^2+28\eqb\eqc^2\eqd+\eqc^4-12\eqc\eqd^2 .
\end{align*}
\end{defi}
\begin{defi}\label{defe-del-equ-bis}
We associate to $\eqpoly$ the following expressions in terms of these given by definition \ref{defe-del-equ}:
\begin{align*}
\deltsl:=  (\eqd-\eqb\eqc)\sqrt{\delo} \left(4\eqb^2\eqc^2-4\eqb\eqc\eqd+2\eqc^3+\eqd^2\right)  +\frac{\sqrt{-3}}{9}\dell  ,
\end{align*}
and
\begin{align*} 
& \arga:= -\frac{2\sqrt{-3}}{3}\left(4\eqb^3\eqc-2\eqd\eqb^2-13\eqb\eqc^2+15\eqd\eqc\right) +2\eqc\sqrt{\delo} ,  \\
& \argb:=   8\eqb^5\eqc^2-8\eqb^4\eqc\eqd -40\eqb^3\eqc^3+2\eqb^3\eqd^2 +116\eqb^2\eqc^2\eqd +23\eqb\eqc^4-99\eqb\eqc\eqd^2 -21\eqc^3\eqd+27\eqd^3    -\sqrt{-3}\left( 8\eqb^2\eqc^2-10\eqb\eqc\eqd+c^3+3\eqd^2\right) \sqrt{\delo}  .  \nonumber
\end{align*}
\end{defi}

Now, we present the principal result of the paper which is the formula of the roots of $\eqpoly$. This formula is in fact a Lagrange-type \cite{Lag}, where the square and cubic are given by the standard convention, contrary to Cardan-Tartaglia and Lagrange formulas which are incorrect with the standard convention. Indeed, the difference with the previous works is the two terms $\prefaca$ and $\prefacb$ of \eq{sol-eq-c-card} which were missed in Cardan-Tartaglia and Lagrange formulas:
\begin{theorem}\label{teo-centrale}
Let the two cubic roots $\raca$ and $\racb$:
\begin{align*} 
& \raca:= \sqrt[3]{  \frac{\sqrt{-3}}{9} \sqrt{\delo} +\frac{2\eqb^3-9\eqc\eqb+27\eqd}{27} }  ,  \nonumber \\
& \racb:= \sqrt[3]{  \frac{\sqrt{-3}}{9} \sqrt{\delo} -\frac{2\eqb^3-9\eqc\eqb+27\eqd}{27} }   .
\end{align*}
The set of solutions of $\eqpoly$ is given in terms of $\raca$ and $\racb$ as follows:
\begin{align}\label{sol-eq-c-card}
 \eqxm =  \indm \prefaca \raca   +  \indm^2 \prefacb  \racb  -\frac{\eqb}{3} ;~ ~~~~
 \indm \in  \left\{ -1 ~,~~  \frac{1-\sqrt{-3}}{2} ~,~~ \frac{1+\sqrt{-3}}{2}\right\} ,
\end{align}
where $\prefaca$ and $\prefacb$ are defined as follows:
\begin{align*} 
& \prefaca : = \frac{\sqrt[3]{4}}{2} \exp\left\{\sqrt{-1}\left[ \arg\left(\arga\sqrt[3]{\deltsl} \right) - \arg\left(-\denumo \raca \right) \right]  \right\} , \\
& \prefacb : = \frac{\sqrt[3]{4}}{2} \exp\left\{\sqrt{-1}\left[ \arg\left(\argb\sqrt[3]{\deltsl^2} \right) - \arg\left( \denumo^2 \racb \right) \right]  \right\}.
\end{align*}
Here, $\delo$, $\denumo$, $\deltsl$, $\arga$ and $\argb$ are given by definitions \ref{defe-del-equ} and \ref{defe-del-equ-bis}.
\end{theorem}
\begin{preuve}
Firstly, we construct the following matrix 
\begin{align*} 
\mata:=\begin{pmatrix}
   -\eqb & 1 & \eqd^{-1} \eqc \\
   0 & 0 & 1 \\
   -\eqd & 0 & 0 
\end{pmatrix},
\end{align*}
which has $\eqpoly$ as its characteristic polynomial. Consequently, the roots of $\eqpoly$ are identical to the spectrum of $\mata$, which can be calculated from corollaries \ref{cor2-p} and \ref{cor3-p}. But, these corollaries require firstly to make explicitly $\paro$ and $\rl$ relative to matrix $\mata$. Indeed, lemma \ref{lem-1} gives \eq{equ-o-teo} of $\paro$ and the expression of $\parl$ in terms of the components of $\matm$. So, we identify the components of $\matm$ with these of $\mata$ so that $\matm=\mata$. In other words, we replace in expressions of section \ref{sec-1} the components of $\matm$ by their corresponding values in terms of $\eqb$, $\eqc$ and $\eqd$. Thus, we deduce that 
\begin{align} \label{def-o-equ}
\paro= \frac{ \eqb^2\eqc^3-\eqc^4-3\eqc\eqd^2+\eqb^2\eqd^2-2\eqb^3\eqc\eqd \pm\sqrt{-3}\eqc(\eqd-\eqb\eqc)\sqrt{\delo} }{\denumo} ,
\end{align} 
since \eq{equ-o-teo} is a quadratic polynomial equation for $\paro$. In the sequel, we take the positive sign before $\sqrt{-3}\eqc(\eqd-\eqb\eqc)\sqrt{\delo}$ in \eq{def-o-equ}. Also, the applying of lemma \ref{lem-1} on matrix $\mata$ implies that
\begin{align*} 
\rl = \lo-\indm\sqrt[3]{4} \frac{\eqc }{\denumo} \sqrt{\delo} \sqrt[3]{\deltsl} ;~ ~~
 \indm \in  \left\{ -1 ~,~~  \frac{1-\sqrt{-3}}{2} ~,~~ \frac{1+\sqrt{-3}}{2}\right\} .
\end{align*}
Here, $\lo$ is defined as follows:
\begin{align*} 
\lo:= \denumo^{-1} \left[ 2\eqb^4\eqc\eqd+\eqb^3\eqc^3-\eqb^3\eqd^2-5\eqb^2\eqc^2\eqd-4\eqb\eqc^4+6\eqb\eqc\eqd^2+5\eqc^3\eqd +\frac{\sqrt{-3}}{3}\eqc\left(\eqb^2\eqc-2\eqb\eqd-\eqc^2\right)\sqrt{\delo} \right]  .           
\end{align*}
Secondly, since we have the expressions of $\paro$ and $\rl$ relative to matrix $\mata$, then by applying corollaries \ref{cor2-p} and \ref{cor3-p} on matrix $\mata$, we deduce that the roots of $\eqpoly$ are given by
\begin{multline}\label{sol-cor2-equ}
\eqxm = \frac{1}{2}\left[\left( 1+\frac{1}{\paro} \right) 
\frac{-\eqc\rl^2 +\eqd\left( 1-\paro\right)\rl +\eqb\eqd\paro\left( 1-\paro\right)+\eqc \left( \eqb^2-4\eqc\right)  \paro^2 }
{\left( 1-\paro\right)\left[ \eqd+(\eqb\eqc-\eqd)\paro\right] -2\eqb\eqc\paro^2 -\eqc\left( 1+\paro\right)\rl}
 -\eqb -\frac{\rl}{\paro}\right] =
 \\
 \frac{1}{2}\left[\left( 1+\frac{1}{\paro} \right) 
 \frac{\eqb\eqd(\eqb\eqc-\eqd) +\eqc\eqd\left(\eqb^2-2\eqc\right)-\eqd\left( 2\eqc^2-\eqb\eqd\right)\paro^{-1} +\eqd\paro^{-1}\left(\eqd\paro^{-1}  -\eqd+2\eqb\eqc\right) \rl}
 {2\eqb^2\eqc^2-4\eqb\eqc\eqd-2\eqc^3+\eqd^2 -2\left( \eqc^3-\eqb\eqc\eqd+\eqd^2\right) \paro^{-1} +\eqd^2\paro^{-2}-2\eqc(\eqd-\eqb\eqc)\paro^{-1}\rl}
  -\eqb -\frac{\rl}{\paro}\right] .
\end{multline} 
Then, we aim to simplify the expression of $\eqxm$ in \eq{sol-cor2-equ}, indeed \eq{sol-cor2-equ} itself implies that 
\begin{multline}\label{sol-eq-c-double}
\frac{-\eqc\rl^2 +\eqd\left( 1-\paro\right)\rl +\eqb\eqd\paro\left( 1-\paro\right)+\eqc \left( \eqb^2-4\eqc\right)  \paro^2 }
{\left( 1-\paro\right)\left[ \eqd+(\eqb\eqc-\eqd)\paro\right] -2\eqb\eqc\paro^2 -\eqc\left( 1+\paro\right) \rl} = \\
\frac{\eqb\eqd(\eqb\eqc-\eqd) +\eqc\eqd\left(\eqb^2-2\eqc\right)-\eqd\left( 2\eqc^2-\eqb\eqd\right)\paro^{-1} +\eqd\paro^{-1}\left(\eqd\paro^{-1}  -\eqd+2\eqb\eqc\right) \rl}
{2\eqb^2\eqc^2-4\eqb\eqc\eqd-2\eqc^3+\eqd^2 -2\left( \eqc^3-\eqb\eqc\eqd+\eqd^2\right) \paro^{-1} +\eqd^2\paro^{-2}-2\eqc(\eqd-\eqb\eqc)\paro^{-1} \rl}.
\end{multline}
So, the following property allows us to eliminate, in the denominator of the fractions of \eq{sol-eq-c-double}, the cubic root concealed in $\rl$:
\begin{align}\label{prop-eq-c-double}
\frac{\ac}{\bc}=\frac{\cc}{\dc}\Rightarrow \frac{\ac}{\bc}=\frac{\cc}{\dc}=\frac{\mc \ac +\nc \cc}{\mc \bc +\nc \dc};~\forall \ac,\bc,\cc,\dc,\mc,\nc\in\esc^*.
\end{align}
By applying property (\ref{prop-eq-c-double}) on \eq{sol-eq-c-double} for $\mc=2\paro^{-1}\eqc(\eqd-\eqb\eqc)$ and $\nc=-\eqc\left( 1+\paro\right)$, we just obtain a square root in the denominator as follows:
\begin{align}\label{prop-eq-c-double-cons}
&\eqxm = -\frac{\eqb}{3}  + \frac{ \indm^2 \sqrt[3]{4^2\deltsl^2} } {  12 \deltsl   } 
\left[\frac{\sqrt{-3}}{3}\left(4\eqb^3\eqc-2\eqd\eqb^2-13\eqb\eqc^2+15\eqd\eqc\right) +\eqc\sqrt{\delo} \right]
- \indm\sqrt[3]{4\deltsl}  \times
\\
&\frac{  8\eqb^5\eqc^2-8\eqb^4\eqc\eqd -40\eqb^3\eqc^3+2\eqb^3\eqd^2 +116\eqb^2\eqc^2\eqd +23\eqb\eqc^4-99\eqb\eqc\eqd^2 -21\eqc^3\eqd+27\eqd^3   
 + \sqrt{-3}\left( 8\eqb^2\eqc^2-10\eqb\eqc\eqd+c^3+3\eqd^2\right) \sqrt{\delo}     }
  {   18  \deltsl  }  .  \nonumber
\end{align}
Therefore, we rationalize \eq{prop-eq-c-double-cons} using the following identities:
\begin{multline*}
 8\eqb^5\eqc^2-8\eqb^4\eqc\eqd -40\eqb^3\eqc^3+2\eqb^3\eqd^2 +116\eqb^2\eqc^2\eqd +23\eqb\eqc^4-99\eqb\eqc\eqd^2 -21\eqc^3\eqd+27\eqd^3   
\\
+ \sqrt{-3}\left( 8\eqb^2\eqc^2-10\eqb\eqc\eqd+c^3+3\eqd^2\right) \sqrt{\delo} =
\frac{9\arga\deltsl}{4\denumo}    , 
\\
\frac{2\sqrt{-3}}{3}\left(4\eqb^3\eqc-2\eqd\eqb^2-13\eqb\eqc^2+15\eqd\eqc\right) +2\eqc\sqrt{\delo} = \frac{3\argb\deltsl}{\denumo^2} 
\end{multline*}
and consequently we deduce that 
\begin{align}\label{sol-eq-cc}
\eqxm = -\frac{\sqrt[3]{4}}{8} \frac{ \indm\denumo \arga \sqrt[3]{\deltsl} - \indm^2\sqrt[3]{4}\argb \sqrt[3]{\deltsl^2} }{ \denumo^2}  -\frac{\eqb}{3} .
\end{align}
Thirdly, in order to simplify $\denumo^2$ in the denominator of \eq{sol-eq-cc}, we prove that:
\begin{align} \label{norme-cub-sol-eq-c}
\begin{split}
 \left(\arga\sqrt[3]{\deltsl} \right)^3 = \left(-4 \denumo \raca \right)^3 \Rightarrow &
  \arga\sqrt[3]{\deltsl}  = \exp\left\{\sqrt{-1}\left[ \arg\left(\arga\sqrt[3]{\deltsl} \right) - \arg\left(-4 \denumo \raca \right) \right]  \right\}  \left(-4 \denumo \raca \right) ,  \\
 \left(\argb\sqrt[3]{\deltsl^2} \right)^3 =  \left(\sqrt[3]{4^2} \denumo^2 \racb \right)^3 \Rightarrow &
  \argb\sqrt[3]{\deltsl^2}  =    \exp\left\{\sqrt{-1}\left[ \arg\left(\argb\sqrt[3]{\deltsl^2} \right) - \arg\left(\sqrt[3]{4^2} \denumo^2 \racb \right) \right]  \right\} \sqrt[3]{4^2} \denumo^2 \racb . 
\end{split}  
\end{align}
Then, by inserting \eq{norme-cub-sol-eq-c} in \eq{sol-eq-cc}, we deduce \eq{sol-eq-c-card}.
\end{preuve}

\section*{Acknowledgments} 
I would like to thank, in particularly, Dalia Ibrahim for her support during the difficult moment to work this paper. Also, I wish to thank the reviewers of the first submission for their useful comments.

\section*{}

\bibliographystyle{elsarticle-num}
\bibliography{biblio}

\end{document}